\documentclass{amsart}

\usepackage{amsmath,amstext,amsfonts,amssymb}
\usepackage{xcolor}

\newtheorem{theorem}{Theorem}
\newtheorem{lemma}{Lemma}
\newtheorem{corollary}{Corollary}

\newcommand{\RR}{\mathbb R}
\newcommand{\NN}{\mathbb N}

\newcommand{\Mm}{\mathcal{M}_{m}} 
\newcommand{\Mw}{\mathfrak{M}_{m}} 

\DeclareMathOperator{\M}{{\bf E}}

\allowdisplaybreaks
\begin{document}

\title{On discretization of some extremal problems}

\author{Oleg Kovalenko}
\email{olegkovalenko90@gmail.com}
\address{Oles Honchar Dnipro National University, Ukraine}

\begin{abstract}
 We solve two continuous extremal problems on the classes of monotone functions: in the first problem we find extremal values for a line integral of a coordinate-wise monotone function of two variables from a rearrange\-ment-invariant class of functions;  in the second one we find extremal values for the expectation of a random process with monotone trajectories at a random time. In both cases we reduce the continuous problems to their discrete counterparts. The obtained discrete problems are on the one hand interesting on their own, and on the other hand give a natural explanation of the structure of the extremal functions for the continuous problems.

\end{abstract}

\maketitle

\section{Introduction}

Recall (see e.g.~\cite{exactConstants}[Section~3.2.2]) that for a measurable function $f\colon E\to\RR$ defined on a measurable space $(E,\mu)$ with a finite measure $\mu$, the function
$$
m_f\colon [0,\infty)\to [0,\mu E],\, m_f(t) = \mu\{x\in E\colon |f(x)| > t\}
$$
is called the distribution function of $f$.
The function
$$
 r_f\colon [0,\mu E]\to[0,\infty),\,r_f(t) = \inf\{s\colon m_f(s)\leq t\} 
$$
is called the non-increasing rearrangement of $f$.

As usual, for two partially ordered sets $(A,\preceq)$ and $(X,\leq)$ a function $f\colon A\to X$ is called non-decreasing if 
$
\alpha\preceq \beta\implies f(\alpha) \leq f(\beta).
$
A function $f\colon [0,1]^2\to\RR$ is called coordinate-wise non-decreasing, if it is non-decreasing with respect to the partial order 
\begin{equation}\label{xyOrder}
    (x_1,y_1)\preceq (x_2,y_2)\iff x_1\leq x_2\text{ and } y_1\leq y_2.
\end{equation}

Let $m\colon [0,1]\to[0,1]$ be an  increasing bijection. Denote by $\Mm$ the class of all measurable (with respect to the two-dimensional Lebesgue measure $\mu$) coordinate-wise non-decreasing  functions 
$f\colon [0,1]^2\to [0,1]$ such that
    \begin{equation}\label{distributionFunction}
        m_f(t) \geq 1-m(t), t\in [0,1].
    \end{equation}
\begin{theorem}\label{th::coordinatewiseMonotone}
   Let $m\colon [0,1]\to[0,1]$ be an  increasing bijection.  For each continuous non-decreasing function $t\colon  [0,1]\to [0,1]$ one has
\begin{equation}\label{linearIntegralEstimate}
    \inf_{f\in\Mm}\int_0^1 f(t(s),s)ds =  \int_0^1 m^{-1}(t(s)\cdot s)ds.
\end{equation}
\end{theorem}

Recall that (see e.g.~\cite[211I, 211J]{FremlinV2}) a set $A\in \mathcal{F}$ of a  measure space $\{\Omega,\mathcal{F},P\}$ is called an atom if $P(A)>0$ and for each $B\in \mathcal{F}$, $B\subset A$, either $P(A\setminus B) =0$ or  $P(B) = 0$. A measure space without atoms is called atomless. 

A natural analogue of the class $\Mm$ for functions of two variables that are monotone with respect to the partial order 
\begin{equation}\label{xMonotonicity}
   (x_1,y_1)\preceq (x_2,y_2)\iff x_1 \leq x_2 \text{ and } y_1 = y_2, 
\end{equation}
can be defined as follows. Assume an atomless probability space $\{\Omega,\mathcal{F},P\}$ is given. For a function $m$ as above, let $\Mw$ be the space of all measurable random processes $\{\xi_t,t\in[0,1]\}$  defined on  $\{\Omega,\mathcal{F},P\}$ such that 
$$
s<t\implies 0\leq \xi_s(\omega)\leq \xi_t(\omega) \leq 1\text{ for almost all } \omega \in \Omega,
$$
and
$$
\M \mu\{t\in [0,1]\colon \xi_t > s\} \geq 1-m(s)\text{ for all } s\in [0,1],
$$
where $\M$ denotes the expectation of a random variable, and $\mu$ is the one-dimensional Lebesgue measure.
\begin{theorem}\label{th::randomProcesses}
Let $m\colon [0,1]\to[0,1]$ be an  increasing bijection, and $\tau$ be a random variable on $\{\Omega,\mathcal{F},P\}$  that takes values in $[0,1]$. Then
\begin{equation}\label{expectationEstimate}
    \inf_{\xi_t\in \Mw} \M \xi_\tau  =  \int_0^1m^{-1}\left(\int_{1-y}^1 r_\tau(s)ds\right)dy.
\end{equation}
 If in addition $m(t) = t$ for all $t\in[0,1]$, then one has
$$
  \inf_{\xi_t\in \Mw} \M\xi_\tau 
  =
\int_0^1 r_\tau(s)s  ds.
$$
\end{theorem}

Classes $\Mm$ and $\Mw$ are examples of rearrangement-invariant classes, which appear in many extremal problems of approximation theory, see e.g.~\cite{Babenko07,Babenko09,Babenko10}. We also note that in the situation, when $\tau$ is not constant, the extremal random process for Theorem~\ref{th::randomProcesses} that we construct in the proof can not be represented as a product 
\begin{equation}\label{tensorProcess}
  \xi_t(\omega) = x(t)\cdot \eta(\omega)  
\end{equation}
 of a function $x\colon [0,1]\to\RR$ and a random variable $\eta$, cf. remarks after~\cite[Theorem~2]{Kovalenko24c}. Solution to some other extremal problems of approximation theory for random processes can be found in~\cite{Drozhzhina,Kovalenko20a,Kovalenko21a}; in all these three articles the extremal processes are of form~\eqref{tensorProcess}.

The proofs of the estimates from below for the quantities  $\int_0^1 f(t(s),s)ds$, $f\in\Mm$, and $\M \xi_\tau$, $\xi_t\in \Mw$ from Theorems~\ref{th::coordinatewiseMonotone} and~\ref{th::randomProcesses} can be obtained directly in the continuous setting; moreover, the extremal functions will also be given explicitly. However, from our point of view, discretization of these extremal problems gives a natural approach to their solutions, explains the origin of the extremal functions, and could probably be used for a solution of some other extremal problems on these or related classes.

In order to introduce the related discrete extremal problems, we illustrate the discretization approach for Theorem~\ref{th::coordinatewiseMonotone} in the partial case when $m(s) = s$, $s\in [0,1]$ and $t(s) = \alpha$, $s\in[0,1]$ for some fixed number $\alpha\in [0,1]$. For large $n\in\NN$  divide the unit square $[0,1]^2$ into $n^2$ equal squares, and approximate a function $f\in \Mm$ by a coordinate-wise non-decreasing piecewise-constant function $f_n\colon [0,1]^2\to\Xi_n$, where
\begin{equation}\label{Xi_n}
    \Xi_n = \left\{\frac 1{n^2}, \frac 2{n^2},\ldots, \frac{n^2}{n^2}\right\}.
\end{equation}
The function $f_n$ can be identified with a bijection from $\{(i,j)\colon 1\leq i,j\leq n\}$ to $\Xi_n$. If $k = k(n)$ is the number of the column that contains the set $\{\alpha\}\times [0,1]$, then $\frac 1n \sum_{\nu=1}^n f_n\left(k,  \nu\right)$ is an approximation of the integral $\int_0^1 f(\alpha,s)ds$. 

Thus we arrive at a discrete extremal problem, which we formulate and study in a somewhat more general situation then it is needed for Theorems~\ref{th::coordinatewiseMonotone} and~\ref{th::randomProcesses}.

Let $N\in\NN$, $A = \{\alpha_1,\ldots, \alpha_N\}$ be  a partially ordered set with a partial order $\preceq$. As usual, two elements $\alpha,\beta\in A$ are called comparable, if either $\alpha\preceq\beta$, or $\beta\preceq\alpha$; otherwise the elements $\alpha,\beta$ are called incomparable.

Let $\Xi$ be a linearly ordered set (i.e., every two elements in $\Xi$ are comparable) with an order $\leq$, and with a commutative associative binary operation $+$ that agrees with the ordering in the sense that 
\begin{equation}\label{+<agreement}
 \xi,\eta,\zeta\in \Xi, \xi <  \eta\implies \xi + \zeta < \eta + \zeta.   
\end{equation}

Assume that $\Xi_N = \{\xi_1,\ldots, \xi_N\}\subset \Xi$ is some subset of $N$ different elements $\xi_1<\xi_2<\ldots <\xi_N$. 
We consider the set $F$ of non-decreasing bijections $f\colon A\to \Xi_N$ and solve the following problem. Given a fixed set of 'indices' $B = \{\beta_1, \ldots, \beta_n\}\subset A$, $n\leq N$, find  quantities
\begin{equation}\label{minMax}
  \min_{f\in F} S(f) \text{ 
 and } \max_{f\in F}S(f),
\end{equation}
where 
$$
S(f) = S(B,f):=\sum_{k=1}^n f(\beta_k).
$$

We need the following notations. 
For a finite set $C$ denote by $|C|$ the number of its elements. Denote by $P_k$ the set of permutations of the set $\{1,\ldots, k\}$, $k\in\NN$ i.e., the set of all bijections $\pi\colon \{1,\ldots, k\}\to \{1,\ldots, k\}$. For each $\alpha \in A$ we set 
$$
\Pi_{\alpha} = \{\beta\in A\colon \beta\preceq \alpha\}
\text { and } 
\Pi^{\alpha} = \{\beta\in A\colon \alpha\preceq \beta\}.
$$
For brevity we also set
$$
\Pi_k = \Pi_{\beta_k} \text { and } \Pi^k = \Pi^{\beta_k}, k=1,\ldots, n.$$ For each $\pi\in P_n$ we define sets 
$$
T_{\pi,k} := \bigcup_{j = 1}^k \Pi_{\pi(j)} \text{ and } T^{\pi,k} := \bigcup_{j = 1}^k \Pi^{\pi(j)}, k = 1,\ldots, n.
$$
Observe that for each $\pi\in P_n$, $f\in F$ and  $k\in \{1,\ldots, n\}$ one has 
\begin{equation}\label{maxOnT}
\max_{T_{\pi,k}} f = \max\{f(\beta_{\pi(1)}), \ldots, f(\beta_{\pi(k)})\}.
\end{equation}
Indeed, if $\alpha\in T_{\pi,k}$, then there exists $j\in \{\pi(1),\ldots, \pi(k)\}$ such that $\alpha\in \Pi_{\pi(j)}$, hence $\alpha\preceq \beta_{\pi(j)}$, and thus $f(\alpha)\leq f(\beta_{\pi(j)})$.

Finally, for each $\pi\in P_n$ such that 
\begin{equation}\label{piRestriction}
1\leq k,j\leq n, \beta_k\prec \beta_j\implies \pi(k) < \pi(j)
\end{equation}
we denote by 
$$
F_\pi = \{f\in F\colon f(\beta_{\pi(1)}) < \ldots <  f(\beta_{\pi(n)}) \}.
$$
In the proof of Theorem~\ref{th::min} below we show that each of these sets $F_\pi$ is non-empty.

Using the introduced notations, we can formulate the following result.
\begin{theorem}\label{th::min}
For each $\pi\in P_n$ that satisfies~\eqref{piRestriction}
\begin{equation}\label{conditionalMin}
\min_{f\in F_\pi }S(B,f) = \sum_{k=1}^n \xi_{\left|T_{\pi,k}\right|}
\text{ and }
\max_{f\in F_\pi }S(B,f) = \sum_{k=1}^n \xi_{N-\left|T^{\pi,k}\right| + 1}
\end{equation}
If $f_*$ and $f^*$ are functions on which the minimum and maximum in~\eqref{conditionalMin} are attained respectively, then  for $k = 1,\ldots, n$
\begin{equation}\label{f(beta)}
f_*\left(\beta_{\pi(k)}\right) = \xi_{|T_{\pi,k}|}\text{ and } 
f^*\left(\beta_{\pi(k)}\right) = \xi_{N - |T^{\pi,k}| + 1}.
\end{equation}
Hence,
\begin{equation}\label{thMin}
\min_{f\in F}S(B,f) = \min_\pi\sum_{k=1}^n \xi_{\left|T_{\pi,k}\right|} \text{ and } 
\max_{f\in F}S(B,f) = \max_\pi\sum_{k=1}^n \xi_{N - |T^{\pi,k}| + 1},
\end{equation}
where the minimum and the maximum on the right-hand sides of~\eqref{thMin} are taken over all $\pi\in P_n$ that satisfy~\eqref{piRestriction}.
\end{theorem}
In the general case, when the set $B$ contains many incomparable elements, the minimum and the maximum on the right-hand sides of~\eqref{thMin} are taken over a too large set of permutations $\pi$; thus direct application of~\eqref{thMin} becomes computationally infeasible. However, under some additional assumptions Theorem~\ref{th::min} may imply explicit results.

\begin{corollary}\label{th::monotoneNodes}
If $\beta_k\prec \beta_{k+1}$ for all $k=1,\ldots, n-1$, then
$$
\min_{f\in F}S(B,f) = \sum_{k=1}^n \xi_{\left|\Pi_{k}\right|} \text{ and } 
\max_{f\in F} S(B,f) = \sum_{k=1}^n \xi_{N-\left|\Pi^{k}\right|+ 1}.
$$
\end{corollary}
\begin{corollary}\label{th::emptyIntervalsIntersection}
If $
|\Pi_1|\leq |\Pi_2|\leq\ldots\leq |\Pi_n|,$
and for all $k\neq j$, $\Pi_k\cap\Pi_j =\emptyset$, then
$$
\min_{f\in F}S(B,f) = \sum_{k=1}^n \xi_{\sum_{j=1}^k\left|\Pi_{j}\right|}.
$$
If $
|\Pi^1|\leq |\Pi^2|\leq\ldots\leq |\Pi^n|,
$
and for all $k\neq j$, $\Pi^k\cap\Pi^j =\emptyset$, then
$$
\max_{f\in F}S(B,f) = \sum_{k=1}^n \xi_{N -\sum_{j=1}^k\left|\Pi^j\right| + 1}.
$$
\end{corollary}

The article is organized as follows. In Section~\ref{s::discreteProblems} we prove the formulated results for the discrete problem. Section~\ref{s::continuousProblems} is devoted to solutions of the continuous problems.

\section{Discrete extremal problems}\label{s::discreteProblems}
\subsection{A connection between the maximization and the minimization problems}\label{s::connection}
The problems of finding minimum in~\eqref{minMax} can in fact be reduced to finding the maximum in~\eqref{minMax} and vice versa. Indeed, if we consider reversed orders in both sets $A$ and $\Xi$ (i.e., introduce an order $\preceq_1$ in the set $A$ such that $\alpha\preceq_1 \beta\iff \beta\preceq \alpha$ and similarly $\leq_1$ in $\Xi$), then a function $f$ is non-decreasing with respect to the initial orders if and only if it is non-decreasing with respect to the new orders. Moreover, $\min$ with respect to the order $\leq$ is $\max$ with respect to the order $\leq_1$. Thus applying results regarding minimum with respect to the initial orders, we obtain results regarding maximum with respect to the new orders and vice versa.

This means that it is sufficient to prove the formulated results only regarding minimum. Below we prove only results for the minimization problem.
\subsection{Auxiliary results}
In order to prove Theorem~\ref{th::min}, we need several auxiliary results.
\begin{lemma}\label{l::twoValuesSubstitution}
Let $f\in F$ and two incomparable elements $\alpha, \beta\in A$  be such that $f(\alpha) = \xi_i$, $f(\beta) = \xi_{i+1}$ for some $i\in\{1,\ldots, N-1\}$. Then the function 
$$
f_{\alpha,\beta}(\gamma) =
\begin{cases}
f(\gamma), & \gamma\notin\{\alpha,\beta\}\\
f(\beta), &\gamma = \alpha\\
f(\alpha),& \gamma = \beta
\end{cases}
$$
also belongs to $F$.
\end{lemma}
\begin{proof}
Assume the contrary, let $\varepsilon,\delta\in A$ be such that 
\begin{equation}\label{nonMonotonicity}
\varepsilon\prec \delta,\text{ but } f_{\alpha,\beta}(\varepsilon) > f_{\alpha,\beta}(\delta).
\end{equation}
 Since $\alpha$ and $\beta$ are incomparable, $\{\alpha,\beta\}\neq \{\varepsilon,\delta\}$. Since the functions $f$ and $f_{\alpha,\beta}$ differ only at two points and $f\in F$, assumption~\eqref{nonMonotonicity} implies $\{\alpha,\beta\}\cap \{\varepsilon,\delta\}\neq \emptyset$. Thus the sets $\{\alpha,\beta\}$ and $\{\varepsilon,\delta\}$ have precisely one common element. We consider the case $\alpha\in \{\varepsilon, \delta\}$, a contradiction in the other case can be obtained analogously. 
 
 If $\alpha = \varepsilon$, then $f(\delta)> f(\varepsilon) = f(\alpha) = \xi_i$. Since $\delta\neq \beta$, and $f$ is a bijection, we obtain that $f(\delta)\neq f(\beta) = \xi_{i+1}$, and hence
 $$
 f_{\alpha,\beta}(\delta) 
 = f(\delta)\geq \xi_{i+2}
 > \xi_{i+1} = f_{\alpha,\beta}(\alpha) 
 =
 f_{\alpha,\beta}(\varepsilon),
 $$
 which contradicts to~\eqref{nonMonotonicity}.

 If $\alpha = \delta$, then 
 $$
 f_{\alpha,\beta}(\varepsilon) = f(\varepsilon) < f(\delta) 
 = f(\alpha) = \xi_i < \xi_{i+1} = f_{\alpha,\beta}(\alpha)
 = f_{\alpha,\beta}(\delta),
 $$
 which again contradicts to~\eqref{nonMonotonicity}. The lemma is proved.
\end{proof}

\begin{lemma}\label{l::extFuncValues}
Let $f$ be a function on which minimum in~\eqref{conditionalMin} is attained.
Then for each $k = 1,\ldots,n$,
\begin{equation}\label{smallestGoFirst}
f\left(T_{\pi,k}\right) = \left\{\xi_1, \ldots, \xi_{|T_{\pi,k}|}\right\}.
\end{equation}
\end{lemma}
\begin{proof}
Assume the contrary, and let $k^*$ be some value of $k\in \{1,\ldots, n\}$ for which equality~\eqref{smallestGoFirst} does not hold. Then there exist $\alpha\in A\setminus T_{\pi,k^*}$, $\beta\in T_{\pi,k^*}$, and $i\in \{1,\ldots, N-1\}$ such that $f(\alpha) = \xi_i$, and $f(\beta) = \xi_{i+1}$.  

We show that 
\begin{equation}\label{alphaIsNotChosen}
  \alpha\notin B.  
\end{equation}
Indeed, if $k^* = n$, then $T_{\pi,k}\supset B$ and~\eqref{alphaIsNotChosen} follows. Otherwise, due to the choice of $\pi$ (which guarantees $f(\beta_{\pi(1)}) < \ldots <  f(\beta_{\pi(n)})$), and   property~\eqref{maxOnT},
\begin{multline*}
 \min f\left(B\setminus T_{\pi,k^*}\right) > \max f\left(B\cap T_{\pi,k^*}\right) = \max f\left(T_{\pi,k^*}\right)
 \\ \geq 
 f(\beta) = \xi_{i+1} > \xi_{i} = f(\alpha),   
\end{multline*}
thus $\alpha\notin B\setminus T_{\pi,k^*}$, and hence~\eqref{alphaIsNotChosen} it true.

Since $f(\alpha)< f(\beta)$, the inequality $\beta\preceq\alpha$ does not hold. Moreover, the inequality $\alpha\preceq\beta$ also can not hold, since otherwise $\alpha$ would belong to $T_{\pi,k^*}$. Thus, the elements $\alpha$ and $\beta$ are incomparable.

Consider the function $f_1 = f_{\alpha,\beta}$ defined in Lemma~\ref{l::twoValuesSubstitution}. By its definition, $f_1(\gamma)\leq f(\gamma)$ for all $\gamma\neq \alpha$. 
Due to~\eqref{alphaIsNotChosen}, we obtain that $f_1(\gamma)\leq f(\gamma)$ for all $\gamma\in B$. Since the minimum in~\eqref{minMax} is attained on the function $f$, due to~\eqref{+<agreement} we obtain $f_1(\gamma) = f(\gamma)$ for all $\gamma\in B.$
Moreover, due to~\eqref{+<agreement}, for the quantity $S(T_{\pi,k^*},f) = \sum_{\gamma\in T_{\pi,k^*}}f(\gamma)$, we obtain $S(T_{\pi,k^*},f_1) < S(T_{\pi,k^*},f)$.

Applying the same arguments to the function $f_1$ we  build $f_2$, then  $f_3$, and so on. For these functions $f_s$, $s\geq 1$ we have
\begin{equation}\label{fsAreOptimal}
f_s(\gamma) = f(\gamma) \text{ for all } \gamma\in B
\end{equation}
and $S(T_{\pi,k^*},f_s) < S(T_{\pi,k^*},f_{s-1})$ for all $s\geq 2$.
This sequence of functions can not be infinite, since the set $\{S(T_{\pi,k^*},g), g\in F\}$ is finite. Thus for some $s\in\NN$ we can not build the function $f_{s+1}$, which means that
$$
f_s\left(T_{\pi,k^*}\right) = \left\{\xi_1, \ldots, \xi_{|T_{\pi,k^*}|}\right\}.
$$
However, due to property~\eqref{maxOnT}, this means that 
$$
\max_{\gamma\in B\cap T_{\pi,k^*}}f_s(\gamma) =
\max_{\gamma\in T_{\pi,k^*}}f_s(\gamma) 
< 
\max_{\gamma\in T_{\pi,k^*}}f(\gamma)
=
\max_{\gamma\in B\cap T_{\pi,k^*}}f(\gamma),
$$
which contradicts to~\eqref{fsAreOptimal}.
\end{proof}
\subsection{Proof of Theorem~\ref{th::min}}\label{thMinProof}
\begin{proof}
Due to Lemma~\ref{l::extFuncValues}, in order to prove~\eqref{conditionalMin} and~\eqref{f(beta)} it is enough to show that $F_\pi$ is non-empty. Indeed, if $F_\pi$ is non-empty, then it contains a finite number of elements, thus the minimum from~\eqref{conditionalMin} exists. Properties~\eqref{smallestGoFirst} and~\eqref{maxOnT} imply~\eqref{f(beta)}, which in turn implies equality~\eqref{conditionalMin}.

We set $T_{\pi, 0} = \emptyset, T_{\pi,n+1} = A$, $A_k = T_{\pi,k}\setminus T_{\pi,k-1}$, $N_k = |T_{\pi,k}|$, $k = 0,\ldots, n+1$. Then the sets $A_k$, $k =1,\ldots, n+1$, are pairwise disjoint, non-empty (except possibly $A_{n+1}$), and $\bigcup_{k=1}^{n+1}A_k = A$; if $A_{n+1}$ is empty, then in the arguments below we exclude the value $k=n+1$. 

For each $k = 1,\ldots, n+1$ consider a bijective non-decreasing function 
$$
f_k\colon A_k\to \{\xi_{N_{k-1} + 1}, \xi_{N_{k-1} + 2},\ldots, \xi_{N_{k}}\},
$$
and define 
$$
f(\gamma) = f_k(\gamma), \gamma \in A_k, k = 1,\ldots, n+1.
$$
Next we show that
\begin{equation}\label{fIsMonotone}
    f\in F.
\end{equation}
 Let $\alpha\prec \beta$ and $\alpha\in A_i, \beta\in A_j$, $i,j\in \{1,\ldots, n+1\}$. 
Then $\beta\in T_{\pi,j}$, and hence due to construction of the set $T_{\pi,j}$, we obtain $\alpha\in T_{\pi,j}$, thus $i\leq j$. If $i = j$, then
$$
f(\beta) = f_j(\beta)\geq f_j(\alpha) = f(\alpha),
$$
due to monotonicity of the function $f_j$. If $i < j$, then 
$$
f(\alpha) = f_i(\alpha) \leq \xi_{N_{i}} < \xi_{N_{j-1}+1}\leq f_j(\beta) = f(\beta),
$$
which finishes the proof of~\eqref{fIsMonotone}.
Thus we obtain that property~\eqref{smallestGoFirst} holds, and hence, due to property~\eqref{maxOnT}, equalities~\eqref{f(beta)} hold. Thus $f\in F_\pi$.
\end{proof}
\subsection{Proof of Corollary~\ref{th::monotoneNodes}}
\begin{proof}
Under conditions of the theorem, minimum in~\eqref{thMin} is taken over a single permutation $\pi$ (which is the identical mapping). Moreover, in this case $T_{\pi,k} = \Pi_k$ for all $k=1,\ldots, n$, so the statement of the theorem follows from~\eqref{thMin}.
\end{proof}
\subsection{Proof of Corollary~\ref{th::emptyIntervalsIntersection}}
\begin{proof}
Under conditions of the theorem the elements of the set $B$ are pairwise incomparable, since $\beta_i\prec \beta_j\implies \Pi_i\subset \Pi_j$. Hence condition~\eqref{piRestriction} holds for arbitrary $\pi\in
 P_n$. Moreover, for arbitrary $\pi\in 
 P_n$, $|T_{\pi,k}|  = \sum_{j=1}^k\left|\Pi_{\pi(j)}\right|$, and hence
 $$
\sum_{k=1}^n \xi_{\left|T_{\pi,k}\right|}
=
\sum_{k=1}^n \xi_{\sum_{j=1}^k\left|\Pi_{\pi(j)}\right|}
\geq 
\sum_{k=1}^n \xi_{\sum_{j=1}^k\left|\Pi_{j}\right|}.
 $$
In order to finish the proof, it is sufficient to observe that the last inequality becomes equality if $\pi$ is the identity permutation.
\end{proof}
\subsection{Some partial cases}
Let $n\in\NN$, $A = \{(i,j)\colon 1\leq i,j\leq n\}$, $\Xi = \Xi_n$ be defined in~\eqref{Xi_n}, and the partial order in $A$ be defined by~\eqref{xyOrder}. Then from Corollary~\ref{th::monotoneNodes} we obtain for each $f\in F$ and $1\leq s\leq n$
\begin{equation}\label{coordinanteWiseDiscreteCase}
\sum_{v = 1}^nf(s,v) \geq \frac 1{n^2}\left(s + 2s + \ldots + ns\right) = \frac{s(n+1)}{2n}.
\end{equation}

Let  $A$ be as above, the order be defined by~\eqref{xMonotonicity}, and
\begin{equation}\label{XiM} 
\Xi:=\left\{m^{-1}\left(\frac{1}{n^2}\right),m^{-1}\left(\frac{2}{n^2}\right),\ldots, m^{-1}\left(\frac{n^2}{n^2}\right)\right\},
\end{equation} 
where $m\colon[0,1]\to[0,1]$ is an increasing bijection. 
If $1\leq s_1\leq\ldots\leq s_n\leq n$ is a fixed set of indices, then for all $f\in F$, due to Corollary~\ref{th::emptyIntervalsIntersection} 
\begin{equation}\label{xMonotoneDiscreteCase}
\sum_{v = 1}^n f(v, s_v)\geq 
\sum_{\nu = 1}^n m^{-1}\left(\frac{s_1 + \ldots + s_\nu}{n^2}\right).
\end{equation}
\section{Continuous extremal problems}\label{s::continuousProblems}
\subsection{On discretization of Theorem~\ref{th::coordinatewiseMonotone}}\label{s::th1Discretization}
In this section we show how the continuous extremal problem from Theorem~\ref{th::coordinatewiseMonotone} can leverage the solution  to the discrete problem from Corollary~\ref{th::monotoneNodes}. A formal proof of Theorem~\ref{th::coordinatewiseMonotone} 
 will be given in the next section. The primary goal here is to give an insight into the structure of an extremal function for the infimum on the left-hand side of~\eqref{linearIntegralEstimate}. In order to somewhat simplify the notations we consider only the partial case, when  $t(s) = \alpha\in (0,1]$ for all $s\in [0,1]$ and
\begin{equation}\label{mIsIdentical}
   m(t) = t.
\end{equation}

Let $f\in \Mm$, $k\in \NN$, and $\varepsilon > 0$ be fixed. Below it is convenient to work with the semi-open square $(0,1]^2$ instead of $[0,1]^2$.
Set
$$
L_0 = \emptyset, L_s = \left\{ (x,y)\in (0,1]^2\colon f(x,y) \leq \frac s k\right\}, s = 1,\ldots, k.
$$
For each $s = 1,\ldots, k$, due to the assumption~\eqref{mIsIdentical},
\begin{equation}\label{muLi}
\mu L_s \leq \frac sk.
\end{equation}
 
We choose a number $\varepsilon_0 < \varepsilon$ such that $\beta = \alpha - \varepsilon_0\in \mathbb{Q}$ and $\beta\geq 0$.
For each $i = 1,\ldots, k$, the intersection of the set $L_i$ with the vertical line $x = \beta$ is a connected set, since $f$ is coordinate-wise monotone. 
Let $j_*$ be the smallest index such that $L_{j_*}$ has a non-empty intersection with the line $x = \beta$, $j^*$ be the smallest index for which $L_{j^*}$ contains the point $(\beta,1)$,   and let $l_s = \sup\{y\colon (\beta,y)\in L_s\}$ be the ordinates of the upper ends of of the intersections, $s = j_*, \ldots, j^*$. It might happen that some of the numbers $l_s$ are equal, but for simplicity of the notations we assume that $l_{j_*} < l_{j_*+1}< \ldots < l_{j^*} =1$.

We choose $\varepsilon_s > 0$ in such a way, that  $l_s - \varepsilon_s\in \mathbb{Q}$, 
$ l_s-\varepsilon_s > l_{s-1}$
 for all $s = j_*,\ldots, j^*$, $l_{j_*-1}:= 0$,    and 
\begin{equation}\label{sumOfEpsilons}
\sum_{s=j_*}^{j^*} \varepsilon_s < \varepsilon.
\end{equation}
Due to coordinate-wise monotonicity of the function $f$, 
\begin{equation}\label{rectangleIsInLevelSet}
\{(x,y)\in (0,1]^2\colon x\leq \beta, y\leq l_s - \varepsilon_s\}\subset L_s, s= j_*,\ldots, j^*.  
\end{equation}
There exists an infinite set of natural numbers $n$ such that each of the numbers
\begin{equation}\label{choiceOfN}
\frac{n}{k}, (l_s-\varepsilon_s)n, s = j_*,\ldots, j^*,  \text{ and } \beta n\text{ are natural}.
\end{equation}
 For each such $n$ we split the square $(0,1]^2$ into $n^2$  squares
 $$
 S_{ij} = \left(\frac{i-1}{n}, \frac i n\right]\times \left(\frac{j-1}{n}, \frac j n\right], i,j = 1,\ldots, n,
 $$
 and build a coordinate-wise non-decreasing bijective function 
 $$
 f_n\colon A := \{(i,j)\colon 1\leq i,j\leq n\} \to\Xi_n,$$ where $\Xi_n$ is defined in~\eqref{Xi_n}. 

The construction of the function $f_n$ is done in $k$ steps. On each step $s\in \{1,\ldots, k\}$ we define the function $f_n$ on the set $A_s\subset A$ of indices such that $S_{ij}\subset L_s$ for all $(i,j)\in A_s$, and where $f_n$ is not yet defined; the values of $f_n$ are chosen in such a way that $f_n$ is injective coordinate-wise non-decreasing and $n^2\cdot f_n (i,j)\in \{1,\ldots, |A_s|\}$ for all $(i,j)\in A_s$. Such construction is done similarly to the one used in Section~\ref{thMinProof}.
Set 
$$
 D_s: = \bigcup_{(i,j)\in A_s} S_{ij}, s = 1,\ldots, k
$$
and 
$$
E_s = \{\beta\}\times (l_{s-1}, l_s-\varepsilon_s], s = j_*, \ldots, j^*.
$$
According to the construction, and in view of~\eqref{rectangleIsInLevelSet} and~\eqref{choiceOfN}, 
$$
E_s\subset D_s\subset L_s\text{ and } E_s\cap L_{s-1} = \emptyset\text{ for each } s = j_*, \ldots, j^*.
$$
Hence $f(x,y) \in \left(\frac{s-1}{k}, \frac sk\right]$ for all $(x,y)\in E_s$, $s = j_*, \ldots, j^*$. Moreover, due to~\eqref{muLi}, for each $s=1,\ldots, k$, $|A_s|\leq \frac sk n^2$, and hence $f_n(i,j)\leq \frac s k$ for all $(i,j)\in A_s$. Thus for all $s = j_*, \ldots, j^*$, $(x,y)\in E_s$, and all $(i,j)\in A_s$ we have 
\begin{equation}\label{fnIsSmaller}
    f(x,y)\geq f_n(i,j) - \frac 1k.
\end{equation}

Define $\varphi_n(x,y) = f_n(i,j)$, $(x,y)\in S_{ij}$, set $\Delta = \bigcup_{s = j_*}^
{j^*}(l_s-\varepsilon_s, l_s]$, and  $\Delta_s := \left(\frac{s-1}n, \frac sn\right]$, $s=1,\ldots, n$. Then the function $\varphi_n(\beta, \cdot)$ is measurable and takes a finite number of values. Moreover, due to~\eqref{fnIsSmaller} and~\eqref{sumOfEpsilons}, inequality
$$
f(\beta,y)\geq \varphi_n(\beta,y) - \frac 1k
$$
holds on all intervals $\Delta_s$ except possibly at most $(2k+ \varepsilon n)$ ones (this is an upper bound for the number of intervals $\Delta_s$ that may intersect the set $\Delta$)  i.e., it holds on a set of (one-dimensional) measure at least $1-\frac{2k}{n}-\varepsilon$.

Thus using monotonicity of $f$,  the fact $f(x,y)\leq 1$ for all $(x,y)\in (0,1]^2$, and~\eqref{coordinanteWiseDiscreteCase},  we obtain
\begin{multline*}
\int_0^1f(\alpha, y)dy+ \frac 1k + \frac{2k}{n} +\varepsilon
\geq
\int_0^1f(\beta, y)dy + \frac 1k + \frac{2k}{n} +\varepsilon
\geq 
\int_0^1\varphi_n(\beta, y)dy
    \\ =
    \frac 1n\sum_{s=1}^nf_n(\beta n, s) 
     \geq
    \frac 1n \cdot \frac{\beta n(n+1)}{2n} 
    \geq   
    \frac{(\alpha-\varepsilon) (n+1)}{2n} 
\end{multline*}
Due to arbitrariness of $\varepsilon>0$ and $n,k\in\NN$,  this implies that 
\begin{equation}\label{th1PartialCase}
    \int_0^1f(\alpha, y)dy\geq \frac \alpha 2,
\end{equation} which is precisely the necessary lower bound in Theorem~\ref{th::coordinatewiseMonotone} in the considered partial case.

The function $n^2\cdot \varphi_n$ can be defined as follows: set it to $1,2,\ldots, \beta n$ on the sets $S_{11},\ldots, S_{(\beta n) 1}$ respectively, then  'fill' the next 'row' of the sets $S_{12},\ldots, S_{(\beta n) 2}$ with the next consecutive natural numbers $\beta n +1, \ldots, 2\beta n$, and so on, until the values on $(0,\beta n]\times (0,1]$ are assigned. Then repeat the same procedure on the set $(\beta n,1]\times (0,1]$, using the numbers $\beta n^2+1, \ldots, n^2$. This implies that
$$
n^2\cdot \varphi_n(x,y) = 
\begin{cases}
[yn]\beta n + O(n)    ,& x\leq \beta\\
\beta n^2 + [yn](1-\beta) n + O(n)    ,& x > \beta
\end{cases}
$$ 
as $n\to\infty$, where $[\gamma]$ denotes the integer part of a real number $\gamma$, and suggests to consider the following function
$$
f(x,y) = \begin{cases}
    \alpha y, & x\leq \alpha\\
    (1-\alpha)y + \alpha,& x > \alpha.
\end{cases}
$$
Direct computations show that $f\in \Mm$ and that $\int_0^1f(\alpha,y) dy = \frac \alpha 2$, which proves that inequality~\eqref{th1PartialCase} is sharp on $\Mm$.

Observe that inequality~\eqref{distributionFunction} was only used to obtain~\eqref{mIsIdentical}, which in turn was used to prove~\eqref{fnIsSmaller}. 
In view of the discussion from Section~\ref{s::connection} it is easy to see that on the class of measurable coordinate-wise non-decreasing  functions 
$f\colon [0,1]^2\to [0,1]$ such that
$m_f(t) \leq 1-t, t\in [0,1]$, using the same arguments one can compute
$$
 \sup_{f}\int_0^1 f(\alpha,y)dy.
$$
The extremal function will be $$
f(x,y) = \begin{cases}
    \alpha y, & x< \alpha\\
    (1-\alpha)y + \alpha,& x \geq \alpha.
\end{cases}
$$
\subsection{Proof of Theorem~\ref{th::coordinatewiseMonotone}}
In this section we give a proof of Theorem~\ref{th::coordinatewiseMonotone}. Observe that both the estimate from below and formula~\eqref{th1Extremal} for an extremal function are essentially continuous versions of the approach outlined in the previous section.
\begin{proof}
    Fix an arbitrary $s\in [0,1]$. Let $f\in \Mm$ and $u = f(t(s),s)$. Then for all $x\leq t(s)$, $y\leq s$ one has $f(x,y)\leq u$, and hence
    $$
    1-m(u)\leq \mu\{(x,y)\in [0,1]^2\colon f(x,y) > u\}\leq 1 - t(s)s.
    $$
    Hence $f(t(s),s) = u\geq m^{-1}(t(s)\cdot s)$, which implies 
\begin{equation}\label{continuousProofEstimateFromBelow}
      \int_0^1 f(t(s),s)ds \geq  \int_0^1 m^{-1}(t(s)\cdot s)ds.  
\end{equation}
   Next we construct a function $f\in \Mm$ for which equality in the latter inequality is attained. For $s\in [0,1]$ set  
    $$\Pi_s = [0, t(s)]\times [0,s]\setminus\bigcup_{u < s} [0, t(u)]\times [0,u].$$
Then $s_1\neq s_2\implies \Pi_{s_1}\cap \Pi_{s_2} = \emptyset$,  $\bigcup_{s\in [0,1]}\Pi_s = [0,t(1)]\times [0,1]$, and we can define a function
\begin{equation}\label{th1Extremal}
f(x,y) = \begin{cases}
    m^{-1}(t(s)\cdot s), & (x,y)\in \Pi_s,\\
    m^{-1}(t(1) + (1-t(1))y), & x > t(1).
\end{cases}
\end{equation}
Since $(t(s), s)\in \Pi_s$ for all $s\in [0,1]$, it is easy to see that inequality~\eqref{continuousProofEstimateFromBelow} becomes equality for function~\eqref{th1Extremal}. Thus to finish the proof it now suffices to show that $f\in \Mm$.

We prove that $f$ is coordinate-wise monotone first.

If $x_1\leq x_2\leq t(1)$, $y_1\leq y_2$, and $(x_2,y_2)\in \Pi_{s_2}$, $s_2\in [0,1]$, then for some $0\leq s_1\leq s_2$ one has $(x_1,y_1)\in \Pi_{s_1}$, and hence $f(x_1,y_1)\leq f(x_2,y_2)$. Thus $f$ is coordinate-wise monotone on $ [0,t(1)]\times [0,1]$. Moreover, it is easy to see that $f$ is coordinate-wise monotone on $(t(1),1]\times [0,1]$. Finally, if $x\leq t(1)$, then $f(x,y)\leq m^{-1}(t(1))$, and $f(x,y)\geq m^{-1}(t(1))$ for all $x > t(1)$. Hence $f$ is coordinate-wise monotone on $[0,1]^2$.

Next we prove that inequality~\eqref{distributionFunction} holds for function~\eqref{th1Extremal}. If $0\leq u\leq t(1)$, then from continuity of the function $s\mapsto st(s)$ it follows that there is a number $s\in [0,1]$ such that $st(s) = u$. Then $f(x,y)\leq m^{-1}(u)$ for all $(x,y)\in \bigcup_{u\leq s}\Pi_u = [0,t(s)]\times [0,s]$ and $f(x,y) > u$ for all other $(x,y)\in [0,1]^2$. Thus 
$$\mu\{(x,y)\colon f(x,y) > m^{-1}(u)\} = 1- st(s) = 1- u.$$
 If $1\geq u > t(1)$, then 
\begin{multline*}
\mu\{(x,y)\colon f(x,y) > m^{-1}(u)\} = 
\mu\{(x,y)\colon x> t(1), t(1) + (1-t(1))y >u\} 
\\= (1-t(1))\cdot \left(1-\frac{u-t(1)}{1-t(1)}\right) = 1-u,
\end{multline*}
as desired.
\end{proof}

\subsection{Proof of Theorem~\ref{th::randomProcesses}}
The estimate from below for Theorem~\ref{th::randomProcesses} can also be proved directly in the continuous situation. However, it is somewhat more complicated then for Theorem~\ref{th::coordinatewiseMonotone}, and the idea of the proof still originates from  the 
 discretization. So we give the proof for the estimate from below using the discretization ideas. Sharpness of the obtained estimate is shown by a random process constructed explicitly in continuous situation; its origin again comes from the discrete considerations.

We use ideas similar to the ones from Section~\ref{s::th1Discretization}, so we omit some details in the proof below.

\begin{proof}
Let $\xi_t\in\Mw$ and $\delta,\varepsilon > 0$ be fixed.
First we approximate the random variable $\tau$ by a simple random variable $\tau_\delta$ in the following sense. There are $v\in\NN$, sets
$\Omega_1,\ldots, \Omega_v\in \mathcal{F}$ such that 
 $\tau_\delta(\omega)  =\tau_i\in [0,1]\cap \mathbb{Q}$ for all $\omega\in \Omega_i$, $i =1,\ldots, v$, $\tau_1\leq \tau_2\leq\ldots\leq \tau_v$, and
$$
\tau_\delta(\omega)\leq \tau(\omega) < \tau_\delta(\omega) + \delta\text { for all } \omega\in \Omega.
$$
For a fixed $k\in\NN$ consider the sets
$$
L_i = \left\{(t,\omega)\colon \xi_t(\omega) \leq m^{-1}\left(\frac ik\right)\right\}, i = 1,\ldots,k.
$$ 
Then 
\begin{equation}\label{muPLi}
\mu\times P(L_i)\leq \frac ik, i = 1,\ldots,k.
\end{equation}
Next we find sets $\Omega_{ij}\subset \Omega_i\cap L_j$, $\Omega_{ij}\in\mathcal{F}$, $P(\Omega_{ij})\in\mathbb{Q}$, $i=1,\ldots, v$, $j = 1,\ldots, k$ such that for $\Omega_0 :=\Omega\setminus\bigcup_{i,j}\Omega_{ij}$,
$P(\Omega_0) < \varepsilon.$
In order to simplify further notations we assume that $\Omega_0 = \emptyset$.

There exists an infinite set of numbers $n\in\NN$ such that all the numbers $\frac{n}{k}$, $\tau_i n$,  and $P(\Omega_{ij}) n$, $i=1,\ldots, v$, $j=1,\ldots, k$ are natural. For each such $n$, we divide each of the sets $\Omega_{ij}$ and the segment $[0,1]$ into 'pieces' with measure $\frac 1n$, and obtain a 'grid' $S_{ij}$ on the set $[0,1]\times \Omega$. This can be done due to atomlessness of $\{\Omega,\mathcal{F},P\}$, see e.g.~\cite[215D]{FremlinV2}. 

Next we build a random process $\zeta_n$ with monotone trajectories and values in the set~\eqref{XiM} similarly to how the function $\varphi_n$ was build in Section~\ref{s::th1Discretization}, 
so that for all $\omega\in\Omega$
\begin{multline*}
\xi(\tau_\delta(\omega),\omega)\geq \zeta_n(\tau_\delta(\omega),\omega) -\max_{s=1,\ldots, k} (m^{-1}(s/k) - m^{-1}((s-1)/k))
\\= \zeta_n(\tau_\delta(\omega),\omega) + o(1), k\to\infty,
\end{multline*}
where the inequality follows from~\eqref{muPLi} similarly to how~\eqref{muLi} was used in Section~\ref{s::th1Discretization}, and the asymptotic equality follows from the uniform continuity of $m^{-1}$ on $[0,1]$.

We set 
$$
s_i = \begin{cases}
     \tau_1n, & i =  1,\ldots,  nP(\Omega_{1}),\\
     \tau_2n, & i =  nP(\Omega_{1}) + 1,\ldots,  nP(\Omega_{1}\cup \Omega_2),\\
     \cdots\\
      \tau_vn, & i = n P\left(\bigcup_{i=1}^{\nu-1} \Omega_i\right) +1, \ldots, n.
\end{cases}
$$
Observe that due to such definition of the numbers $s_i$,  one can easily prove by induction on $\nu$ that for each $i = 1,\ldots, n$
$$
\sum_{j=1}^i\frac{s_j}{n^2} = \int_{1-\frac in}^1 r_{\tau_\delta}(s)ds.
$$
One has 
\begin{gather*}
  \M\xi_{\tau} 
  = 
\int_{\Omega}\xi(\tau(\omega),\omega)P(d\omega) 
  \geq \int_{\Omega}\xi(\tau_\delta(\omega),\omega)P(d\omega)
   \geq 
\int_{\Omega}\zeta_n(\tau_\delta(\omega),\omega)P(d\omega) + o(1)
 \\\stackrel{\eqref{xMonotoneDiscreteCase} }{\geq}
 \frac 1n\sum_{i=1}^nm^{-1}\left(\sum_{j = 1}^i \frac{s_j}{n^2}\right)  + o(1)
  =
  \frac 1n\sum_{i=1}^nm^{-1}\left(\int_{1-\frac in}^1 r_{\tau_\delta}(s)ds\right)  + o(1)
   \\
  \stackrel{n\to\infty}{\to}\int_{0}^1m^{-1}\left(\int_{1-y}^1r_{\tau_\delta}(s)ds\right)dy  +o(1)
  \stackrel{\delta\to 0}{\to}
    \int_0^1m^{-1}\left(\int_{1-y}^1 r_\tau(s)ds\right)dy + o(1).
\end{gather*}
This implies an estimate from below in~\eqref{expectationEstimate}.
Next we prove that it is sharp.

We first assume that for almost all $\omega\in\Omega$
\begin{equation}\label{tauIsInjective}
\tau(\omega)\in [0,1)\text { and } P\{\eta\in \Omega\colon \tau(\eta) = \tau(\omega)\} = 0.
\end{equation}
Then the function $r_\tau$ is strictly decreasing, and for all $s\in [0,1]$, 
\begin{equation}\label{P<=r(s)}
P\{\omega\colon \tau(\omega)\leq r_\tau( s)\} = 1-s.
\end{equation}
For each $\omega\in\Omega$ we set $\Pi_\omega = \{\eta\in \Omega\colon \tau(\eta)\leq \tau(\omega)\}$. 
Define 
$$
\xi^*_t(\omega) = 
\begin{cases}
    m^{-1}\left(\int_{1-P(\Pi_\omega)}^1 r_\tau(s)ds\right), & t\leq \tau(\omega),\\
    m^{-1}\left(\M\tau + \int_{1-P(\Pi_\omega)}^1 (1-r_\tau(s))ds\right), & t> \tau(\omega).
\end{cases}
$$
For each $\omega\in \Omega$,  $\xi_t^*(\omega)$ is constant on $[0,\tau(\omega)]$ and on $(\tau(\omega), 1]$, and for $0\leq s\leq \tau(\omega) < t\leq 1$, 
$$
m(\xi^*_s(\omega)) = \int_{1-P(\Pi_\omega)}^1 r_\tau(s)ds 
\leq 
\int_{0}^1 r_\tau(s)ds = \M\tau\leq m(\xi^*_t(\omega));
$$
hence $\xi^*_t$ has non-decreasing trajectories. Moreover, $m\circ \xi_0^*$ attains all values from the segment $[0, \M\tau]$, and $m\circ\xi_1^*$ attains all values from the segment $[\M\tau, 1]$.
If $0\leq s\leq \M\tau$ and $\omega^*$ is such that $\int_{1-P(\Pi_{\omega^*})}^1 r_\tau(u)du = s$, then 
\begin{multline*}
\mu\times P\left(\{(t,\omega)\colon \xi^*_t(\omega)\leq m^{-1}(s)\}\right)
=
\mu\times P\left(\{(t,\omega)\colon \omega\in \Pi_{\omega^*}, t\leq \tau(\omega)\}\right)
\\=
\int_{1-P(\Pi_{\omega^*})}^1 r_\tau(u)du = s.
\end{multline*}
If $1\geq s> \M\tau$ and $\omega^*$ is such that $\int_{1-P(\Pi_{\omega^*})}^1 (1-r_\tau(u))du = s - \M\tau$, then
\begin{gather*}
\mu\times P\left(\{(t,\omega)\colon \xi^*_t(\omega)\leq m^{-1}(s)\}\right)
\\=
\mu\times P\left(\{(t,\omega)\colon t\leq \tau(\omega)\}\cup \{(t,\omega)\colon \omega\in \Pi_{\omega^*}, t > \tau(\omega)\}\right)
\\=
\M\tau + \int_{1-P(\Pi_{\omega^*})}^1 (1-r_\tau(u))du = s.
\end{gather*}
Thus $\xi^*_t\in\Mw$.
Next, we compute the value $\M\xi^*_\tau$. 

Based on $P$ and $\tau$, we define the image measure $Q$ on $[0,1]$ with respect to the mapping $\varphi(\omega) = P(\Pi_\omega)$, see e.g.~\cite[Definition~7.7]{Schilling}. This means that for each measurable set $B\subset [0,1]$,
$
Q(B) = P(\varphi^{-1}(B)).
$
Then, due to~\eqref{P<=r(s)}, for all $t\in [0,1]$, 
$$
    Q([0,t]) = P\{\omega\colon P(\Pi_\omega) \leq t\}
    \\=
    P\{\omega\colon \tau(\omega) \leq r_\tau( 1-t)\} = t.
$$
Hence $Q$ is the Lebesgue measure on $[0,1]$. Using~\cite[Theorem~14.1]{Schilling}, we obtain
$$
\M\xi^*_\tau 
= 
\int_\Omega m^{-1}\left(\int_{1-P(\Pi_\omega)}^1 r_\tau(s)ds\right)P(d\omega)
 = 
 \int_0^1m^{-1}\left(\int_{1-y}^1 r_\tau(s)ds\right)dy.
 $$
 If in addition $m(t) = t$ for all $t\in[0,1]$, then changing the order of integration, one has
$$
  \M\xi^*_\tau 
  =
  \int_0^1\int_{1-y}^1 r_\tau(s)dsdy=
\int_0^1\int_{1-s}^1 r_\tau(s)dy  ds
=
\int_0^1 r_\tau(s)s  ds,
$$
which finishes the proof of the theorem under assumptions~\eqref{tauIsInjective}.

In the general case we can approximate $\tau$ by random variables $\tau_n$, $n\in\NN$, that converge to $\tau$ uniformly on $\Omega$, and such that~\eqref{tauIsInjective} holds for each of them.
\end{proof}

\bibliographystyle{elsarticle-num}
\bibliography{bibliography}
\end{document}